\documentclass[a4paper, 11pt]{article}
\usepackage{amsmath}
\usepackage{amsfonts}
\usepackage{amssymb}
\usepackage[english]{babel}
\usepackage{graphicx}
\usepackage{amsthm}
\usepackage[title]{appendix}

\usepackage{longtable}
\usepackage{tabularx}
\allowdisplaybreaks

\newtheorem{theorem}{Theorem}[section]

\newtheorem{lemma}[theorem]{Lemma}

\newtheorem{construction}[theorem]{Construction}

\newtheorem{problem}[theorem]{Problem}

\def\whitebox{{\hbox{\hskip 1pt
 \vrule height 6pt depth 1.5pt
 \lower 1.5pt\vbox to 7.5pt{\hrule width
    3.2pt\vfill\hrule width 3.2pt}%
 \vrule height 6pt depth 1.5pt
 \hskip 1pt } }}
\def\qed{\ifhmode\allowbreak\else\nobreak\fi\hfill\quad\nobreak
     \whitebox\medbreak}

\newcommand{\ignore}[1]{}

\begin{document}
\baselineskip 16pt
\title{Further results on the Hamilton-Waterloo problem }

\author{\small  L. Wang, S. Lu, and H. Cao \thanks{Research
supported by the National Natural Science Foundation of China
under Grant 11571179, the Natural Science Foundation of Jiangsu
Province under Grant No. BK20131393, and the Priority Academic
Program Development of Jiangsu Higher Education Institutions.
E-mail: {\sf caohaitao@njnu.edu.cn}} \\
\small Institue of Mathematics, \\ \small  Nanjing Normal
University, Nanjing 210023, China}

\date{}
\maketitle
\begin{abstract}

In this paper, we almost completely solve the existence of an
almost resolvable cycle system with odd cycle length. We also use
almost resolvable cycle systems as well as other combinatorial
structures to give some new solutions to the Hamilton-Waterloo
problem.


\smallskip
\end{abstract}

\section{Introduction}

In this paper, the vertex set and the edge set of a graph $H$ will
be denoted by $V(H)$ and $E(H)$, respectively. We denote the cycle of length $k$ by $k$-cycle
 or $C_k$,  the complete graph on $v$ vertices by $K_v$ and
denote the complete $u$-partite graph with $u$ parts of
size $g$ by $K_u[g]$.

A {\it factor} of
a graph $H$ is a spanning subgraph of $H$.
Suppose $G$ is a subgraph of a graph $H$, a {\it $G$-factor} of graph $H$ is a set
 of edge-disjoint subgraphs of $H$,
each isomorphic to $G$. And a {\it $G$-factorization} of $H$ is a set of edge-disjoint $G$-factors of $H$.
 A $C_k$-factorization of a graph $H$ is a partition of
$E(H)$ into $C_k$-factors.  In
\cite{AH,ASSW, HS, JLA, L, LD,PWL, R}, we can obtain
the following result.

\begin{theorem}\label{Kv}
There exists a $C_k$-factorization of $K_u[g]$ if and only if $
g(u-1 )\equiv 0\pmod 2$, $gu\equiv 0\pmod k$, $k$ is even when
$u=2$, and $(k,u,g)\not\in\{(3,3,2),(3,6,2),(3,3,6),(6,2,6)\}$.
\end{theorem}

An {\it $r$-factor} is
a factor which is $r$-regular. It's obvious that a 2-factor
consists of a collection of disjoint cycles.
 A {\it
$2$-factorization} of a graph $H$ is a partition of $E(H)$ into
2-factors.
The well-known {\it Hamilton-Waterloo problem} is the problem of
determining whether $K_v$ (for $v$ odd) or $K_v$ minus a
$1$-factor (for $v$ even) has a $2$-factorization in which there
are exactly $\alpha$ $C_m$-factors and $\beta$ $C_n$-factors.
The authors {\rm \cite{WCC}} generalize this problem to a general graph $H$, and
use HW$(H;m,n;\alpha,\beta)$ to denote a $2$-factorization of $H$
in which there are exactly $\alpha$ $C_m$-factors and $\beta$
$C_n$-factors. So when $H=K_v$(for $v$ odd) or $K_v$ minus a
$1$-factor (for $v$ even), an HW$(H;m,n;\alpha,\beta)$ is a
solution to the original Hamilton-Waterloo problem, denoted by
HW$(v;m,n;\alpha,\beta)$. Further, denote by
HWP$(v;m,n)$ the set of $(\alpha,\beta)$ for which a solution
HW$(v;m,n;\alpha,\beta)$ exists.
It is easy to see that the necessary conditions for the existence
of an HW$(v;m,n;\alpha,\beta)$ are $m | v$ when $ \alpha >0$, $n |
v$ when $ \beta >0$ and $\alpha+\beta=\lfloor \frac{v-1}{2}
\rfloor$. When $\alpha\beta=0$, there is a solution to
HW$(v;m,n;\alpha,\beta)$, see \cite{AH,ASSW,HS}. From now on, we
suppose that $\alpha\beta\not=0$.

A lot of work has been done for small values of $m$ and $n$.
A complete solution for the existence of an
HW$(v;3,n;\alpha,\beta)$ in the cases $n\in\{4,5,7\}$ is given in
\cite{ABBE,DQS,LF,OO,WCC}.
For the case
$(m,n)\in\{(3,15),(5,15),(4,6),(4,8),(4,16),(8,16)\}$, see
\cite{ABBE}. The existence of an HW$(v;4,n;\alpha,\beta)$ for odd
$n\geq3 $ has been solved except possibly when $v=8n$ and
$\alpha=2$, see \cite{KO,OO,WCC}. It is shown in \cite{K} that the
necessary conditions for the existence of an
HW$(v;3,9;\alpha,\beta)$ are also sufficient except possibly when
$\beta=1$. Many infinite classes of HW$(v;3,3x;\alpha,\beta)$s are
constructed in \cite{AKK}. Much attention to the Hamilton-Waterloo
problem has been dedicated to the case of triangle factors and
Hamilton factors, the results for this case can be found in
\cite{DL,DL2,HNR,LS}.
In \cite{LFS}, the authors give a complete solution for the
existence of an HW$(v;4k,v;\alpha,\beta)$.
Recently, Burgess et al. \cite{BDT} have made significant progress
on the Hamilton-Waterloo problem for uniform odd cycle factors.

\begin{theorem}
\label{odd} If $m$, $n$ and $t$ are odd integers with $n\geq m \geq 3$, then $(\alpha,\beta)\in$ \rm{HWP}$ (mnt;m,n)$ if and only if $\alpha, \beta \geq 0$ and $\alpha+\beta=\frac{mnt-1}{2}$,
 except possibly when:

  $\bullet$ $ t>1$ and $\beta=1$ or $3$, or $(m,n,\beta)=(5,9,5),(5,9,7),(7,9,5),(7,9,7),(3,13,5)$;

  $\bullet$ $ t=1$ and $\beta \in [1,2,\cdots,\frac{n-3}{2}]\cup \{\frac{n+1}{2}, \frac{n+5}{2}\},$
  $(m,\alpha)=(3,2),(3,4)$  or $(m,n,\alpha,\beta)=(3,11,6,10),(3,13,8,11),(5,7,9,8),(5,9,11,11),(5,9,13,9),(7,9,20,11)$ or $(7,9,22,9)$.
\end{theorem}

We continue to consider this problem, the final main results will be the following.

\begin{theorem}
\label{odd2} If $m$, $n$ and $t$ are odd integers with $n\geq m \geq 3$, then $(\alpha,\beta)\in$ \rm{HWP}$ (mnt;m,n)$ if and only if $\alpha, \beta \geq 0$ and $\alpha+\beta=\frac{mnt-1}{2}$,
 except possibly when:

  $\bullet$ $ t>1$ and $\beta=1$ or $3$;

  $\bullet$ $ t=1$ and $\beta \in [1,2,\cdots,\frac{n-3}{2}]\cup \{\frac{n+1}{2}, \frac{n+5}{2}\},$
  $(m,\alpha)=(3,2),(3,4)$.
\end{theorem}

When $m=k$ and $n=2kt+1$, we can also get the following results.

\begin{theorem}
\label{ARCS} Let $J=\{4,6\} \cup \  [8, 9, \cdots, k(k-1)t+\frac{k-3}{2}]$.
If $t \geq 1$ and odd integer $k\geq 3$,
then $(\alpha,\beta)\in$ \rm{HWP}$ (k(2kt+1); k, 2kt+1)$ in the following four cases:

$(1)$ $k=3$: $t\neq 1,2$, $\beta \in [4, 5, \cdots, 6t-6] \cup \{6t\}$ is even;

$(2)$ $k=5$: $\beta \in J \setminus \{ 20t-3, 20t-1 \}$;

$(3)$ $k=7,9$: $\beta \in J $;

$(4)$ $k\geq 11$: $t\neq2$, $\beta \in J $.
\end{theorem}

\begin{theorem}
\label{4k-4kt}  If $v\equiv 0\pmod {4kt} $, then $(\alpha,\beta)\in$ \rm{HWP}$ (v;4k,4kt)$
if and only if $ k \geq 1$, $ t \geq 2$, $\alpha, \beta \geq 0$ and $\alpha+\beta=\frac{v-2}{2}$.
\end{theorem}

In Section 2, we will introduce the definition of an almost
resolvable $k$-cycle system, then we give some constructions for
almost resolvable cycle systems with odd orders. In Section 3, we
will use almost resolvable $k$-cycle system as well as other
combinatorial structures to give some new recursive constructions
for HW$(v;m,n;\alpha,\beta)$ with $m\equiv n\equiv 1\pmod 2$. In
the last section, we shall present some direct constructions and
use these recursive constructions in Section 3 to prove our main
results.

\section{Almost resolvable $k$-cycle systems}

A $k$-{\it cycle system} of order $v$ is that a collection of
$k$-cycles which partition the edges of $K_v$. A $k$-cycle system
of order $v$ exists if and only if $ 3 \leq k \leq v$, $v \equiv 1
\pmod 2$ and $v(v-1) \equiv 0 \pmod {2k}$ \cite{AG,SM}.
A $k$-cycle system of order $v$ is {\it resolvable} if it has a $C_k$-factorization.
 By Theorem~\ref{Kv}, a resolvable $k$-cycle system of order $v$
exists if and only if $3 \leq k \leq v$, $v$ and $k$ are odd, and
$v \equiv 0 \pmod k$. If $v \equiv 1 \pmod {2k}$, then the
$k$-cycle system is not resolvable. In this case, Vanstone et al.
\cite{VSS} started the research of the existence of an almost
resolvable $k$-cycle system.

A collection of $(v-1)/k$ disjoint $k$-cycles is called an {\it
almost parallel class}. In a $k$-cycle system of order $v \equiv 1
\pmod {2k}$, the maximum possible number of almost parallel
classes is $(v-1)/2$ in which case a {\it half-parallel class}
 containing $(v-1)/2k$ disjoint $k$-cycles is left over. A
$k$-cycle system of order $v$ whose cycle set can be partitioned
into $(v-1)/2$ almost parallel classes and a half-parallel class
is called an {\it almost resolvable $k$-cycle system}, denoted by
$k$-ARCS$(v)$.

For recursive constructions of almost resolvable $k$-cycle
systems, C. C. Lindner, et al. \cite{LMR} have considered the
general existence problem of almost resolvable $k$-cycle system
from the commutative quasigroup for $k\equiv0\pmod2$ and make a
 hypothesis:
 if there exists a {\rm $k$-ARCS$(2k+1)$} for $k\equiv0\pmod2$
 and $k\geq8$, then  there exists a {\rm
$k$-ARCS$(2kt+1)$} except possibly for $t=2$. H. Cao et al.
\cite{CNT,NC} continue to consider the recursive constructions of
an almost resolvable $k$-cycle system  for $k\equiv1\pmod2$. By
using recursive method and direct constructions, some classes of
almost resolvable cycle systems with small orders have been
obtained. The known results on the existence of an almost
resolvable cycle system of order $n$ are summarized as below.

\begin{theorem} {\rm (\cite{ABHL,BHL,CNT,DLM,DLR,LMR,VSS})}
\label{3-14} Let $k \geq 3$, $t \geq 1$ be integers and $n = 2kt +
1$. There exists a $k$-{\rm ARCS}$(n)$ for $k \in \{3, 4, 5, 6, 7,
8, 9, 10, 14\}$, except for $(k, n) \in \{(3, 7), (3, 13), (4,
9)\}$ and except possibly for $(k, n) \in \{(8, 33), (14, 57)\}$.
\end{theorem}

\begin{theorem} {\rm (\cite{NC})}
\label{11-49} There exists a $k$-{\rm ARCS}$(2kt+1)$ for $t\geq
1$, $11 \leq k \leq49$, $k \equiv 1 \pmod 2$ and $t\neq 2,3,5$.
\end{theorem}

In this section we focus on constructions of almost resolvable
cycle system with odd orders. The main idea is to find some
initial cycles with special properties such that all the required
almost parallel classes can be obtained from them. We also need
the following notions for our constructions.

Suppose $\Gamma$ is an additive group and  $I=\{\infty\}$ is a set
which is disjoint with $\Gamma$. We will consider an action of
$\Gamma$ on $\Gamma \ \cup \ I$ which coincides with the {\it
right regular action} on the elements of $\Gamma$, and the action
of $\Gamma$ on $I$ will coincide with the identity map. Given a
graph $H$ with vertices in $\Gamma\ \cup\ I$, the \emph{translate}
of $H$ by an element $\gamma$ of $\Gamma$ is the graph $H+\gamma$
obtained from $H$ by replacing each vertex $x\in V(H)$ with the
vertex $x+\gamma$. The {\it development} of $H$ under a subgroup
$\Sigma$ of additive group $\Gamma$ is the collection
$dev_\Sigma(H)=\{H+x\;|\;x\in \Sigma\}$ of all translates of $H$
by an element of $\Sigma$. The \emph{list of differences} of a
graph $H$ with vertices in $\Gamma$ is the multiset $\Delta H$ of
all possible differences $x - y$  with $(x,y)$ an ordered pair of
adjacent vertices of $H$.

\subsection{ Constructions for $k$-ARCS$(2k+1)$}

For our constructions, we suppose
$\Gamma=Z_{u}\times Z_2$ and $I=\{\infty\}$.  For a graph $H$ with
vertices in $\Gamma \cup I$ and any pair $(j,j')\in Z_2\times
Z_2$,  we define \emph{list of $(j,j')$-differences of $H$} as the
multiset $\Delta_{(j,j')}$ of all possible differences $x-y\in
Z_{u}$ with $(x,j)$ adjacent to $(y,j')$ in $H$.

\begin{lemma}\label{A}
Let $v=2k+1$ and let $F$ be a vertex-disjoint union of two cycles of length $k$ satisfying the following conditions:\\
$(i)$ $V(F)=(Z_{k}\times Z_{2})\cup\{\infty\}\backslash \{(a,b)\}$, $(a,b)\in Z_{k}\times Z_{2}$;\\
$(ii)$ $\infty$  has a neighbor in $Z_{k}\times\{0\}$ and the other neighbor in $Z_{k}\times\{1\}$;\\
$(iii)$ $\Delta_{(0, 0)}F\supseteq Z_{k}\setminus\{0\}$,
$\Delta_{(0, 1)}F\supseteq Z_{k}$, $\Delta_{(1, 1)}F\supseteq
Z_{k}\setminus\{0,\pm d\}$, $(d,k)=1$.
\\ Then, there exists a $k$-{\rm ARCS}$(2k+1)$.
\end{lemma}

\noindent{\it Proof:\ } Let $V(K_{v})=(Z_{k}\times
Z_{2})\cup\{\infty\}$.  The half parallel class is the cycle
$C_{0}=((0,1),(d,1)$, $(2d,1),\cdots,(kd-d,1))$ since $(d,k)=1$. By $(i)$, we know that $F$ is an
almost parallel class. All the required $k$ almost parallel
classes will be generated from $F$ by $(+1\pmod{k},-)$.

Now we show that the half parallel class and the $k$ almost
parallel classes form a $k$-ARCS$(2k+1)$.
 Let  $F'=\{(a,1),(a+d,1) \ | \  a\in Z_k\}$ and
 $\Sigma:=Z_{k}\times\{0\}$. Let  $\mathcal{F} = dev_{\Sigma}(F)\ \cup \ \{F'\}$. The
total number of edges-counted with their respective
multiplicities-covered by the almost parallel classes and half parallel class of $\cal F$ is
$k(2k+1)$, that is exactly the size of $E(K_v)$. Therefore, we only
need to prove that every pair of vertices  lies in a suitable
translate of $F$ or in $F'$. By $(ii)$,  an edge $\{(z, j),
\infty\}$ of ${K}_v$ must appear in a cycle of $ dev_{\Sigma}(F)$.

Now consider an edge $\{(z, j), (z', j')\}$ of ${K}_v$ whose
vertices both belong to $Z_{k}\times Z_{2}$. If $j=j'=1$ and $z-z'
\in \{ \pm d\}$, then this edge belongs to $F'$. In all other cases there
is,  by $(iii)$, an edge of $F$ of the form $\{(w, j), (w', j')\}$
such that $w-w'=z-z'$. It then follows that $F+(-w'+z', 0)$ is an
almost parallel class of $ dev_{\Sigma}(F)$ containing the edge
$\{(z, j), (z', j')\}$ and the assertion follows.\qed

Here we use Lemma~\ref{A} to give a construction of a
$k$-ARCS$(2k+1)$ for any  odd $ k\geq 9$.

\begin{lemma}
\label{2k+1} For any odd $ k\geq 9$, there exists a $k$-{\rm ARCS}$(2k+1)$.
\end{lemma}

\noindent {\it Proof:}
Let the vertex set be $(Z_{k} \times Z_{2})\cup \{\infty\}$
and $(a,b)=(0,1)$. For the element $d$ and
two initial cycles $C_1$ and $C_2$ in $F$, we distinguish in two
cases $k\equiv 1\pmod 4$ and $k\equiv 3\pmod 4$.

\vspace{5pt}
\noindent $(1)$ $k\equiv 1\pmod 4$. Let $k=4n+1$ and $d=2n$.

The cycle $C_1$ is the concatenation of the sequences $T_1$, $(0,0)$ and $T_2$ defined
as follows:\\
{\footnotesize
$T_1=((n,0), ( -n,1), \cdots,\underline{(n-i,0),(-(n-i),1)}, \cdots,(1,0),(-1,1))$;\\
$T_2=((1,1),(-1,0),\cdots,\underline{(1+i,1),(-(1+i),0)},\cdots,(n,1),(-n,0))$.}

The cycle $C_2$ is the concatenation of the sequences $\infty$, $T_1$ and $T_2$, where \\
{\footnotesize
$T_1=((-2n,0),  (2n,0),\cdots,\underline{(-(2n-i),0),  (2n-i,0)}, \cdots,(-(n+1),0),(n+1,0))$;\\
$T_2=((n+1,1),(-(n+1),1),\cdots,\underline{(n+1+i,1),(-(n+1+i),1)},\cdots,(2n,1),(-2n,1))$.}

\vspace{5pt}
\noindent $(2)$  $k\equiv 3\pmod 4$. Let $k=4n+3$ and $d=2$.

The cycle $C_1$ is the concatenation of the sequences $T_1$, $T_2$ and $T_3$, where \\
{\footnotesize
$T_1=((-2,0),(2,1),\cdots,\underline{(-(2+2i),0),(2+2i,1)},\cdots,(-2n,0), (2n,1))$;\\
$T_2=((2n,0), (-2n,1),\cdots,\underline{(2n-2i,0),(-(2n-2i),1)},\cdots,(2,0),(-2,1))$; \\
$T_3=((-1,0),(1,1),(0,0))$.}

The cycle $C_2$ is the concatenation of the sequences $\infty$, $T_1$, $T_2$ and $T_3$ defined as
follows: \\
{\footnotesize
$T_1=((2n+1,0),(-(2n+1),0),\cdots,\underline{(2n+1-2i,0),(-(2n+1-2i),0)},\cdots,(3,0),(-3,0))$;\\
$T_2=((1,0),(-1,1))$;\\
$T_3=((3,1),(-3,1),\cdots,\underline{(3+2i,1),(-(3+2i),1)},\cdots,(2n+1,1),(-(2n+1),1))$.}

It is straightforward to check that $F$ satisfies all the conditions of
Lemma~\ref{A}. For brevity, we don't list the vertices and
differences in $F$. The proof is complete.
\qed

\subsection{ Constructions for $k$-ARCS$(6k+1)$}

We construct a $k$-ARCS$(6k+1)$ for any odd $ k \geq 9$.
We give the following lemma which is similar to Lemma~\ref{A}.

\begin{lemma}\label{B}
Let $v=6k+1$ and let $F$ be a vertex-disjoint union of six cycles of length $k$ satisfying the following conditions:\\
$(i)$ $V(F)=(Z_{3k}\times Z_{2})\cup\{\infty\}\backslash \{(a,b)\}$, $(a,b)\in Z_{3k}\times Z_{2}$;\\
$(ii)$ $\infty$  has a neighbor in $Z_{3k}\times\{0\}$ and the other neighbor in $Z_{3k}\times\{1\}$ ;\\
$(iii)$ $\Delta_{(0, 0)}F\supseteq Z_{3k}\setminus\{0\}$,
$\Delta_{(0, 1)}F\supseteq Z_{3k}$,
$\Delta_{(1, 1)}F\supseteq Z_{3k}\setminus\{0,\pm 3\}$.\\
Then, there exists a $k$-{\rm ARCS}$(6k+1)$.
\end{lemma}

\noindent{\it Proof:\ } Let $V(K_{v})=(Z_{3k}\times
Z_{2})\cup\{\infty\}$.  The half-parallel class is $\{C_0 +(i,0) \
| \ i=0,1,2 \}$, where $C_0=((0,1),(3,1),(6,1),
\cdots,(3(k-1),1))$. $F$ is an almost
parallel class by $(i)$. All the required $3k$ almost parallel classes will
be generated from $F$ by $(+1\pmod{3k},-)$.

Next, we show that the half parallel class and the $3k$ almost
parallel classes form a $k$-ARCS$(6k+1)$.
 Let  $F'=\{(a,1),(a+3,1) \ | \ a\in Z_{3k}\}$ and
 $\Sigma:=Z_{3k}\times\{0\}$. Let  $\mathcal{F} = dev_{\Sigma}(F)\ \cup \ \{F'\}$. The
total number of edges-counted with their respective
multiplicities-covered by the almost parallel classes and half parallel class of $\cal F$ is
$3k(6k+1)$, that is exactly the size of $E(K_v)$. So, we only
need to prove that every pair of vertices  lies in a suitable
translate of $F$ or in $F'$. By $(ii)$,  an edge $\{(z, j),
\infty\}$ of ${K}_v$ must appear in a cycle of $ dev_{\Sigma}(F)$.

Next, we consider an edge $\{(z, j), (z', j')\}$ of ${K}_v$ whose
vertices both belong to $Z_{3k}\times Z_{2}$. If $j=j'=1$ and
$z-z' \in \{\pm 3\}$, this edge belongs to $F'$. In all other cases
there is,  by $(iii)$, an edge of $F$ of the form $\{(w, j), (w',
j')\}$ such that $w-w'=z-z'$. It then follows that $F+(-w'+z', 0)$
is an almost parallel class of $ dev_{\Sigma}(F)$ containing the
edge $\{(z, j), (z', j')\}$. Hence, the assertion holds.\qed

\begin{lemma}
\label{6k+1-1} For any $ k\geq 9$ and $k \equiv 1 \pmod 4$, there exists a $k$-{\rm ARCS}$(6k+1)$.
\end{lemma}

\noindent {\it Proof:} Let the vertex set be $(Z_{3k} \times
Z_{2}) \cup \{\infty\}$ and $(a,b)=(0,1)$. Let $k=4n+1$, $n\ge 2$.
The six initial cycles in $F$ are listed as below.

First of all, we construct the following cycles $C_i$, $1 \leq i \leq 4$, where\\
{\footnotesize $C_1=((3n+1,0),(-(3n+1),0),\cdots,\underline{(3n+1+i,0),(-(3n+1+i),0)},\cdots,(5n,0),(-5n,0), (5n+3,0))$; \\
$C_2=((3n+1,1),(-(3n+1),1),\cdots,\underline{(3n+1+i,1),(-(3n+1+i),1)},\cdots,
(5n,1),(-5n,1),(5n+3,1))$;\\
$C_3=((n+1,0),(-(n+1),1),\cdots,\underline{(n+1+i,0),(-(n+1+i),1)},\cdots,
(3n,0),(-3n,1), (-n,0))$;\\
$C_4=((n+1,1),(-(n+1),0),\cdots,\underline{(n+1+i,1),(-(n+1+i),0)},\cdots,
(3n,1),(-3n,0),(-n,1))$.}

\vspace{5pt}

To construct $C_5$ and $C_6$, we start with $k \in \{9,13,17\}$:

{\footnotesize \vspace{5pt}
\noindent
\begin{tabular}{lllll}
$k=9$:
&$C_5=(\infty,(2,1),(1,1),(-1,1),(1,0),(-13,1),(0,0),(-1,0),(2,0))$;\\
&$C_6=((-13,0),(-11,0),(11,1),(12,0),(-12,1),(-12,0),(12,1),(11,0),(-11,1))$.\\
$k=13$:
&$C_5=(\infty,(3,1),(-2,1),(2,1),(1,1),(-1,1),(1,0)$,$(-19,1),(0,0),(-1,0),(2,0),(-2,0), (3,0))$;\\
&$C_6=((-19,0),(-17,0),(17,1),(16,0),(-16,1),(-18,0)$,$(-18,1),(17,0),(-17,1),(19,0)$,\\
& \ \ \ \ \ \ \ \ $(16,1),(-16,0),(19,1))$.\\
$k=17$:
&$C_5=(\infty,(4,1),(-3,1),(3,1),(-2,1),(2,1),(1,1),(-1,1),(1,0),(-25,1),(0,0),(-1,0),(2,0)$,\\
&  \ \ \ \ \ \ \ \  $(-2,0),(3,0),(-3,0),(4,0))$;\\
&$C_6=((-25,0),(-23,0),(22,1),(-22,0),(24,1),(25,0),(21,1),(-21,0),(-24,1),(-24,0)$, \\
&  \ \ \ \ \ \ \ \  $(-21,1),(21,0),(25,1),(24,0),(-22,1),(22,0),(-23,1))$.
\end{tabular}}

\vspace{5pt}

For $k>17$, the cycle $C_5$ is the concatenation of the
sequences $\infty$, $T_1$, $T_2$ and $T_3$, where\\
{\footnotesize
$T_1=((n,1),(-(n-1),1),(n-1,1),\cdots,\underline{(-(n-1-i),1),(n-1-i,1)},\cdots,(-2,1),(2,1)
$,$(1,1),(-1,1))$;\\
$T_2=((1,0),(-(6n+1),1),(0,0),(-1,0))$;\\
$T_3=((2,0),(-2,0),\cdots,\underline{(2+i,0),(-(2+i),0)},\cdots,(n-1,0),(-(n-1),0),(n,0))$.}

\vspace{5pt}

For the cycle $C_6$, we consider the following three cases.

\vspace{5pt}

{\bf Case 1:} $k\equiv1\pmod{12}$,  $ k\geq 25$.

The cycle $C_6$ is the concatenation of the sequences $T_1,
T_2,\ldots, T_{12}$ defined as follows:\\
{\footnotesize $T_1=((-(6n+1),0),(-(6n-1),0))$;\\
$T_2=((6n-5,1),(-(6n-4),0),\cdots,\underline{(6n-5-3i,1),(-(6n-4-3i),0)},\cdots,(5n+4,1)$,$(-(5n+5),0))$;\\
$T_3=((5n+1,1),(-(5n+1),0))$;\\
$T_4=((5n+6,1),(-(5n+4),0),\cdots,\underline{(5n+6+3i,1),(-(5n+4+3i),0)},\cdots,(6n-3,1)$,$(-(6n-5),0))$;\\
$T_5=((6n+1,1),(-(6n-2),0),(6n-1,1),(6n,0),(-6n,1),(6n-2,0))$;\\
$T_6=((-(6n-3),1),(6n-4,0),\cdots,\underline{(-(6n-3-3i),1),(6n-4-3i,0)},\cdots,(-(5n+6),1)$,$(5n+5,0))$;\\
$T_7=((-(5n+3),1),(5n+2,0),(-(5n+2),1),(-(5n+2),0),(5n+2,1),(-(5n+3),0))$;\\
$T_8=((5n+5,1),(-(5n+6),0),\cdots,\underline{(5n+5+3i,1),(-(5n+6+3i),0)},\cdots,(6n-4,1)$,$(-(6n-3),0) )$;\\
$T_9=((6n-2,1),(-6n,0),(6n,1),(6n-1,0),(-(6n-2),1),(6n+1,0))$;\\
$T_{10}=((-(6n-5),1),(6n-3,0),\cdots,\underline{(-(6n-5-3i),1),(6n-3-3i,0)},\cdots,(-(5n+4),1)$,$(5n+6,0))$;\\
$T_{11}=((-(5n+1),1),(5n+1,0))$;\\
$T_{12}=((-(5n+5),1),(5n+4,0),\cdots,\underline{(-(5n+5+3i),1),(5n+4+3i,0)},\cdots,(-(6n-4),1),(6n-5,0)$,
\\ \hspace*{1cm} $(-(6n-1),1))$.}

\vspace{5pt}

{\bf Case 2:} $k\equiv5\pmod{12}$,  $ k\geq 29$.

The cycle $C_6$ is the concatenation of the sequences $T_1,
T_2,\ldots, T_{12}$, where \\
{\footnotesize $T_1=((-(6n+1),0),(-(6n-1),0))$;\\
$T_2=((6n-6,1),(-(6n-5),0),\cdots,\underline{(6n-6-3i,1),(-(6n-5-3i),0)},\cdots,(5n+4,1)$,$(-(5n+5),0))$;\\
$T_3=((5n+1,1),(-(5n+1),0))$;\\
$T_4=((5n+6,1),(-(5n+4),0),\cdots,\underline{(5n+6+3i,1),(-(5n+4+3i),0)},\cdots,(6n-4,1)$,$(-(6n-6),0))$;\\
$T_5=((6n-2,1),(-6n,0),(6n,1),(6n-1,0),(-(6n-3),1),(6n-3,0),(-(6n-2),1)$,$(6n+1,0))$;\\
$T_6=((-(6n-4),1),(6n-5,0),\cdots,\underline{(-(6n-4-3i),1),(6n-5-3i,0)},\cdots,(-(5n+6),1)$,$(5n+5,0))$;\\
$T_7=((-(5n+3),1),(5n+2,0),(-(5n+2),1),(-(5n+2),0),(5n+2,1),(-(5n+3),0))$;\\
$T_8=((5n+5,1),(-(5n+6),0),\cdots,\underline{(5n+5+3i,1),(-(5n+6+3i),0)},\cdots,(6n-5,1)$,$(-(6n-4),0))$;\\
$T_9=((6n+1,1),(-(6n-2),0),(6n-3,1),(-(6n-3),0),(6n-1,1),(6n,0),(-6n,1)$,$(6n-2,0))$;\\
$T_{10}=((-(6n-6),1),(6n-4,0),\cdots,\underline{(-(6n-6-3i),1),(6n-4-3i,0)},\cdots,(-(5n+4),1)$,$(5n+6,0))$;\\
$T_{11}=((-(5n+1),1),(5n+1,0))$;\\
$T_{12}=((-(5n+5),1),(5n+4,0),\cdots,\underline{(-(5n+5+3i),1),(5n+4+3i,0)},\cdots,(-(6n-5),1),(6n-6,0)$,
\\ \hspace*{1cm} $(-(6n-1),1))$.}

\vspace{5pt}

{\bf Case 3:} $k\equiv9\pmod{12}$,  $ k\geq 21$.

The cycle $C_6$ is the concatenation of the sequences $T_1, T_2, \cdots,T_4$,
$(6n+1,0)$, $T_5$, $T_6$, $T_7$, $(6n+1,1)$, $T_8, T_9, \cdots,T_{11}$ defined as follows:\\
{\footnotesize $T_1=((-(6n+1),0),(-(6n-1),0),(6n-1,1),(6n,0),(-6n,1))$;\\
$T_2=((6n-4,0),(-(6n-3),1),\cdots,\underline{(6n-4-3i,0),(-(6n-3-3i),1)},\cdots,(5n+4,0)$,$(-(5n+5),1))$;\\
$T_3=((5n+1,0),(-(5n+1),1))$;\\
$T_4=((5n+6,0),(-(5n+4),1),\cdots,\underline{(5n+6+3i,0),(-(5n+4+3i),1)},\cdots,(6n-2,0)$,$(-(6n-4),1))$;\\
$T_5=((-(6n-2),1),(6n-3,0),\cdots,\underline{(-(6n-2-3i),1),(6n-3-3i,0)},\cdots,(-(5n+6),1)$,$(5n+5,0))$;\\
$T_6=((-(5n+3),1),(5n+2,0),(-(5n+2),1),(-(5n+2),0),(5n+2,1),(-(5n+3),0))$;\\
$T_7=((5n+5,1),(-(5n+6),0),\cdots,\underline{(5n+5+3i,1),(-(5n+6+3i),0)}, \cdots, (6n-3,1)$,$(-(6n-2),0))$;\\
$T_8=((-(6n-4),0),(6n-2,1),\cdots,\underline{(-(6n-4-3i),0),(6n-2-3i,1)},\cdots,(-(5n+4),0)$,$(5n+6,1))$;\\
$T_9=((-(5n+1),0),(5n+1,1))$;\\
$T_{10}=((-(5n+5),0),(5n+4,1),\cdots,\underline{(-(5n+5+3i),0),(5n+4+3i,1)},\cdots,(-(6n-3),0)$,$(6n-4,1))$;\\
$T_{11}=((-6n,0),(6n,1),(6n-1,0),(-(6n-1),1))$.}

We can easily check that $F$ satisfies all the conditions of
Lemma~\ref{B}. \qed

\begin{lemma}
\label{6k+1-3} For any $ k\geq 11$ and $k \equiv 3 \pmod 4$, there exists a $k$-{\rm ARCS}$(6k+1)$.
\end{lemma}

\noindent {\it Proof:}
Let the vertex set be $(Z_{3k} \times Z_{2}) \cup \{\infty\}$ and $(a,b)=(0,1)$. Let $k=4n+3$, $n\ge 2$. The six initial cycles
 in $F$ are given as below.

 We first define the first four cycles $C_i$, $1\leq i \leq4$, where\\
{\footnotesize $C_1=((n+2,0),(-(n+2),0),\cdots,\underline{(n+2+i,0),(-(n+2+i),0)},\cdots,(3n+2,0)$,$(-(3n+2),0),(-(n-1),0))$;\\
$C_2=((n+2,1),(-(n+2),1),\cdots,\underline{(n+2+i,1),(-(n+2+i),1)},\cdots,(3n+2,1)$,$(-(3n+2),1),(-(n-1),1))$;\\
$C_3=((5n+3,0),(-(5n+3),1),\cdots,\underline{(5n+3-i,0),(-(5n+3-i),1)},\cdots,(3n+3,0)(-(3n+3),1)$, \\  \hspace*{1cm} $(-(5n+4),0))$;\\
$C_4=((5n+3,1),(-(5n+3),0),\cdots,\underline{(5n+3-i,1),(-(5n+3-i),0)},\cdots,(3n+3,1),(-(3n+3),0)$,
\\  \hspace*{1cm} $(-(5n+4),1))$.}

To construct $C_5$ and $C_6$, we first deal with $k \in \{11,15,19\}$:

{\footnotesize \vspace{5pt}
\noindent
\begin{tabular}{lllll}
$k=11$:
&$C_5=(\infty,(-2,1),(2,1),(1,1),(3,1),(-3,0),(-3,1),(3,0),(1,0),(-2,0),(2,0))$;\\
&$C_6=((16,0),(-16,0),(15,1),(-15,0),(14,1),(15,0),(-15,1),(14,0),(16,1),(0,0), (-16,1))$.\\
$k=15$:
&$C_5=(\infty,(2,1),(3,1),(-3,1),(1,1),(-1,1),(4,1),(-4,0),(-4,1),(4,0),(1,0),(-1,0),(3,0)$,\\
& \ \ \ \ \ \ \ \  $(-3,0),(2,0))$;\\
&$C_6=((22,0),(-22,0),(21,1),(-21,0),(20,1),(-20,0),(19,1),(20,0),(-20,1),(19,0),(-22,1)$,\\
& \ \ \ \ \ \ \ \  $(0,0),(22,1),(21,0),(-21,1))$.\\
$k=19$:
&$C_5=(\infty,(-2,1),(2,1),(1,1),(3,1),(-4,1),(4,1),(-1,1),(5,1),(-5,0),(-5,1),(5,0),(-2,0)$,\\
& \ \ \ \ \ \ \ \  $(1,0),(3,0),(-1,0),(4,0),(-4,0),(2,0))$;\\
&$C_6=((28,0),(-28,0),(27,1),(-27,0),(26,1),(-26,0),(25,1),(-25,0),(24,1),(25,0),(-25,1)$,\\
& \ \ \ \ \ \ \ \  $(24,0),(-27,1),(27,0),(-26,1),(26,0),(28,1),(0,0),(-28,1))$.
\end{tabular}}

\vspace{5pt}

For $k>19$, we first consider the cycle $C_5$.
It is the concatenation of the sequences $\infty$,
$T_1, T_2, \cdots,T_5$, there are the following three cases.

\vspace{5pt}

{\bf Case 1:} $k \equiv 3\pmod{12}$,  $ k\geq 27$.\\
{\footnotesize $T_1=((2,1),(3,1),(-3,1),(1,1),(-1,1),(4,1))$;\\
$T_2=((-4,1),(6,1),(-6,1),(5,1),(-2,1),(7,1),\cdots,\underline{(-(4+3i),1),(6+3i,1),(-(6+3i),1)(5+3i,1)}$,\\
   \hspace*{1cm}   $\underline{(-(2+3i),1),(7+3i,1)}$,$\cdots,(-(n-2),1),(n,1),(-n,1),(n-1,1)$,$(-(n-4),1),(n+1,1))$;\\
$T_3=((-(n+1),0),(-(n+1),1))$; \\
$T_4=((n+1,0),(-(n-4),0),(n-1,0),(-n,0),(n,0),(-(n-2),0),\cdots,\underline{(n+1-3i,0),(-(n-4-3i),0)}$,\\
  \hspace*{1cm}   $\underline{(n-1-3i,0),(-(n-3i),0),(n-3i,0),(-(n-2-3i),0)},\cdots$,
     $(7,0),(-2,0),(5,0),(-6,0),(6,0),  \\ \hspace*{1cm} (-4,0),(4,0))$; \\
$T_5=((1,0),(-1,0),(3,0),(-3,0),(2,0))$.}

\vspace{5pt}

{\bf Case 2:} $k \equiv 7\pmod{12}$, $ k\geq 31$.\\
{\footnotesize $T_1=((-2,1),(2,1),(1,1),(3,1),(-4,1),(4,1),(-1,1),(5,1))$;\\
$T_2=((-7,1),(7,1),(-3,1),(6,1),(-5,1),(8,1),\cdots,\underline{(-(7+3i),1),(7+3i,1),(-(3+3i),1),(6+3i,1)}$,\\
\hspace*{1cm}
     $\underline{(-(5+3i),1),(8+3i,1)},\cdots,(-n,1),(n,1),(-(n-4),1),(n-1,1)$,
     $(-(n-2),1),(n+1,1))$;\\
$T_3=((-(n+1),0),(-(n+1),1))$;\\
$T_4=((n+1,0),(-(n-2),0),(n-1,0),(-(n-4),0),(n,0),(-n,0),\cdots,\underline{(n+1-3i,0),(-(n-2-3i),0)}$, \\ \hspace*{1cm}
     $\underline{(n-1-3i,0),(-(n-4-3i),0),(n-3i,0),(-(n-3i),0)},\cdots$,
     $(8,0),(-5,0),(6,0),(-3,0),(7,0), \\ \hspace*{1cm} (-7,0),(5,0))$;\\
$T_5=((-1,0),(4,0),(-4,0),(3,0),(1,0),(-2,0),(2,0))$.}

\vspace{5pt}

{\bf Case 3:} $k \equiv 11\pmod{12}$, $ k\geq 23$. \\
{\footnotesize $T_1=((-2,1),(2,1),(1,1),(3,1))$;\\
$T_2=((-5,1),(5,1),(-1,1),(4,1),(-3,1),(6,1),\cdots,\underline{(-(5+3i),1),(5+3i,1),(-(1+3i),1),(4+3i,1)}$, \\
\hspace*{1cm}
     $\underline{(-(3+3i),1),(6+3i,1)},\cdots,(-n,1), (n,1),(-(n-4),1),(n-1,1)$,
     $(-(n-2),1),(n+1,1))$;\\
$T_3=((-(n+1),0),(-(n+1),1))$;\\
$T_4=((n+1,0),(-(n-2),0),(n-1,0),(-(n-4),0),(n,0),(-n,0),\cdots,\underline{(n+1-3i,0),(-(n-2-3i),0),}$
\hspace*{1cm}
     $\underline{(n-1-3i,0),(-(n-4-3i),0),(n-3i,0),(-(n-3i),0)},\cdots$,
     $(6,0),(-3,0),(4,0),(-1,0),(5,0)$,\\ \hspace*{1cm} $(-5,0))$;\\
$T_5=((3,0),(1,0),(-2,0),(2,0))$.}

\vspace{5pt}

$C_6$ is the concatenation of the sequences
$T_1,T_2,\cdots,T_5$. We distinguish two cases.

{\bf Case 1:} $k \equiv 3\pmod8$, $ k\geq 27$.\\
{\footnotesize $T_1=((6n+4,0),(-(6n+4),0))$;\\
$T_2=((6n+3,1),(-(6n+3),0),\cdots,\underline{(6n+3-i,1),(-(6n+3-i),0)},\cdots,(5n+5,1)$,
     $(-(5n+5),0))$;\\
$T_3=((5n+4,1),(5n+5,0),(-(5n+5),1),(5n+4,0))$;\\
$T_4=((-(5n+7),1),(5n+7,0),(-(5n+6),1),(5n+6,0),\cdots,\underline{(-(5n+7+2i),1),(5n+7+2i,0)}$, \\ \hspace*{1cm}
     $\underline{(-(5n+6+2i),1),(5n+6+2i,0)},\cdots,(-(6n+3),1),(6n+3,0)$,
     $(-(6n+2),1), (6n+2,0))$;\\
$T_5=((6n+4,1),(0,0),(-(6n+4),1))$.}

{\bf Case 2:} $k \equiv 7 \pmod 8 $,  $ k \geq 23 $.\\
{\footnotesize $T_1=((6n+4,0),(-(6n+4),0))$;\\
$T_2=((6n+3,1),(-(6n+3),0),\cdots,\underline{(6n+3-i,1),(-(6n+3-i),0)},\cdots,(5n+5,1)$,
     $(-(5n+5),0))$;\\
$T_3=((5n+4,1),(5n+5,0),(-(5n+5),1),(5n+4,0))$;\\
$T_4=((-(5n+7),1),(5n+7,0),(-(5n+6),1),(5n+6,0),\cdots,\underline{(-(5n+7+2i),1),(5n+7+2i,0)}$, \\ \hspace*{1cm}
     $\underline{(-(5n+6+2i),1),(5n+6+2i,0)},\cdots,(-(6n+2),1),(6n+2,0)$,
     $(-(6n+1),1),(6n+1,0))$;\\
$T_5=((-(6n+4),1),(0,0),(6n+4,1),(6n+3,0),(-(6n+3),1))$.}

 $F$ satisfies all the conditions of
Lemma~\ref{B} and the assertion holds.\qed

\subsection{ Constructions for $k$-ARCS$(2kt+1)$}

For the main recursive constructions, we need the definition of
cycle frames. A 2-regular subgraph of a complete multipartite
graph covering all vertices except those belonging to one part $G$
is said to be a {\it holey $2$-factor} missing $G$. We will also
frequently speak of a {\it holey $C_k$-factor} to mean a holey
$2$-factor whose cycles have length $k$.

Let $H$ be a graph $K_u[g]$ with $u$ parts
$G_1,G_2,\ldots,G_u$. A partition of $E(H)$ into holey $2$-factors
of $H-G_i(1\le i\le u)$ is said to be a {\it cycle frame of type
$g^u$}. Further, if all holey 2-factors of a cycle frame of type
$g^u$ are $C_k$-factors, then we denote the cycle frame by
$(k,1)$-cycle frame of type $g^u$.
We write it as $(k, 1)$-CF$(g^u)$ for brevity.
Many authors have contributed to prove the following results.

\begin{theorem}{\rm (\cite{BCDT,BH,CNT,ER, HLR, MS, S, ZC})}
\label{gu}
 There exists a $(k, 1)$-{\rm CF}$(g^u)$ if and only if $ g\equiv 0\pmod 2$, $g(u-1)\equiv 0\pmod
k$, $u\geq 3$ when $k$ is even, $u\geq 4$ when $k$ is odd, except
a $(6,1)$-{\rm CF}$(6^3)$.
\end{theorem}

H. Cao et al. \cite{CNT} give the following general recursive
construction for an almost resolvable cycle system of order $n$
by using cycle frames.

\begin{construction}{\rm (\cite{CNT})}\label{CNT}
Suppose there exists a $(k,1)$-{\rm CF}$(2k)^t$ and a
{\rm $k$-ARCS}$(2k+1)$. Then there exists a {\rm
$k$-ARCS}$(2kt+1)$.
\end{construction}

\begin{theorem}
\label{2kt+1} For any odd $k\geq 11$, there exists a {\rm
$k$-ARCS$(2kt+1)$}, where $t\ge 1$ and $t\neq 2$.
\end{theorem}

\noindent {\bf Proof:} When $t=1,3$, the conclusion comes from
Lemmas~\ref{2k+1},~\ref{6k+1-1} and~\ref{6k+1-3}.

When $t\geq 4$, for any odd $k\geq 5$, there
exists a $(k,1)$-CF$(2k)^t$ by Theorem~\ref{gu}. Applying
Construction~\ref{CNT} with a $k$-ARCS$(2k+1)$ which exists by
Lemma~\ref{2k+1}, we can obtain the required $k$-ARCS$(2kt+1)$.\qed

Further, we use the ARCSs as a tool to obtain some results of the
Hamilton-Waterloo problem. Frist, we give some recursive constructions.

\section{Some new constructions}

Before giving their constructions, we still need
some definitions in graph theory. For more general concepts of
graph theory, see \cite{W}.

Given a graph $G$, $G[n]$ is the
lexicographic product of $G$ with the empty graph on $n$ points.
Specifically, the point set is $\{x_i:x\in V(G), i\in Z_n\}$ and
$(x_i,y_j) \in E(G[n])$ if and only if $(x,y) \in E(G),i,j \in Z_n$. In
the following we will denote by $C_m[n]$ the lexicographic product
of $C_m$ with the empty graph on $n$ points.

We need the following known results for our recursive
constructions.

\begin{theorem}{\rm (\cite{CNT,PWL})}
\label{Cm}  There exists a $C_{m}$-factorization of $C_m[n]$ for $m \geq 3$ and $n \geq 1$
except only for $(m,n)=(3,6)$ and $(m,n) \in \{(l,2) \ | \ l \geq 3 $ is odd $\}$.
\end{theorem}

\begin{theorem}{\rm (\cite{LD})}
\label{Cmn}  There exists a $C_{mn}$-factorization of $C_m[n]$ for $m \geq 3$ and $n \geq 1$.
\end{theorem}

For the next recursive construction, we need more notations.
 When $g(u-1)\equiv 1\pmod 2$, it is easy to see that
 an HW$( K_u[g]; m,n; \alpha,  \beta)$ can not exist.
 In this case, by simple computation, we know that
 it is possible to partition $E(K_u[g])$ into a 1-factor,
 $\alpha$ $C_m$-factors and $\beta$  $C_n$-factors where
 $\alpha+\beta=\lfloor \frac{g(u-1)}{2} \rfloor$.
For brevity, we still use HW$( K_u[g]; m,n; \alpha,  \beta)$
 to denote such a decomposition. Further, denote by
HWP$(K_u[g];m,n)$ the set of $(\alpha,\beta)$ for which a solution
HW$(K_u[g];m,n;\alpha,\beta)$ exists.
Similarly, HWP$(C_u[g];m,n)$ denote the set of $(\alpha,\beta)$ for which a solution
HW$(C_u[g];m,n;\alpha,\beta)$ exists.

\begin{construction} {\rm (\cite{WCC})}
\label{C-RGDD} If there exist an {\rm
HW}$( K_u[g] ; m ,n; \alpha,  \beta)$ and an {\rm
HW}$(g;m,n; \alpha', \beta' )$, then an {\rm
HW}$(gu;m,n;  \alpha + \alpha',  \beta +\beta' )$ exists.
\end{construction}

\begin{construction} {\rm (\cite{WCC})}
\label{L351} If $(\alpha,\beta)\in$\rm{HWP}$(K_u[g];m,n)$,
$(t_i,s-t_i)\in$\rm{HWP}$(C_m[s];m',n')$, $1\le i\le \alpha$, and
$(r_j,s-r_j)\in$\rm{HWP}$(C_n[s];m',n')$, $1\le j\le \beta$, then
$(\alpha',\beta')\in$\rm{HWP}$(K_u[gs];m',n')$, where
$\alpha'=\sum^{\alpha}_{i=1} t_i +\sum^{\beta}_{j=1} r_j$ and
$\beta'=(\alpha+\beta)s-\alpha'$.
\end{construction}

\begin{theorem}{\rm (\cite{BDT})}
\label{Cn-m}  If $m$ and $n$ are odd integers with $n \geq m \geq 3$, $0 \leq \alpha \leq n$,
then $(\alpha,\beta) \in $ {\rm HWP}$(C_m[n]; m, n)$, if and only if $\beta = n-\alpha$, except possibly
when $\alpha=2,4$, $\beta=1,3$, or $(m, n, \alpha)= (3,11,6), (3, 13, 8), (3, 15, 8)$.
\end{theorem}

Let $\Gamma$ be a finite additive group and let $S$ be a subset of
$\Gamma \backslash \{0\}$ such that the opposite of every element
of $S$ also belongs to $S$. The Cayley graph over $\Gamma$ with
connection set $S$, denoted by $Cay(\Gamma, S)$, is the graph with
vertex set $\Gamma$ and edge set $E(Cay(\Gamma, S))= \{(a,b)| a,b
\in \Gamma, a-b \in S\}$. It is quite obvious that $Cay(\Gamma,
S)=Cay(\Gamma, \pm S)$.

\begin{theorem}{\rm (\cite{WCC})}
\label{2p-533} Let $n \geq 3$.  If $a\in Z_n$, the order of $a$
is greater than $3$ and $(i,m)=1$, then there is a $C_m$-factorization of
$Cay(Z_m \times Z_n,  \{ \pm i\}  \times  (\pm \{0, a,  2a\} ) )$.
\end{theorem}

\begin{theorem}{\rm (\cite{BDT})}
\label{d1d2} Let $n \geq m \geq 3$ be odd integers and let $0 < d_1,d_2 < n$.
If any linear combination of $d_1$ and $d_2$ is coprime to $n$, then there exist
four $C_n$-factors which form a $C_n$-factorization of
$Cay(Z_m \times Z_n,  \{ \pm 1\}  \times  (\pm \{ d_1,  d_2\}) )$.
\end{theorem}

\begin{theorem}{\rm (\cite{BDT})}
\label{d}  Let $n \geq m \geq 3$ be odd integers and let $0 < d < n$ be coprime to $n$. There exist
two $C_n$-factors which form a $C_n$-factorization of
$Cay(Z_m \times Z_n,  \{ \pm 1\}  \times  \{ \pm d\} )$.
\end{theorem}

\begin{lemma}
\label{Cm[n]2} $(2,n-2)\in$ \rm{HWP}$(C_m[n];m,n)$ for any odd
integers $n \geq m\geq3$ and $n \equiv 3 \pmod{6} \geq 9$.
\end{lemma}

\noindent {\it Proof:}   Let the vertex set be $\Gamma=Z_{m}\times
Z_n$ and $n=3d$. Let
$C_j=((0,0),(1,b_{j1}),\cdots,(m-1,b_{j,m-1}))$, $1\leq j\leq 2$,
where
$$ b_{11}=-b_{21}=d, \  b_{12}=-b_{22}=2d, \ b_{jt}=b_{j,(t-2)}, \ t\geq 3. $$
Then each $C_j$ will generate a $C_m$-factor by $(-, +1
\pmod{n})$. Thus we get two $C_m$-factors which
form  a $C_m$-factorization of $Cay(\Gamma, \{\pm
1\}\times\{ \pm d\})$.

For $n-2$ $C_n$-factors, five of which can be obtained
from a $C_n$-factorization of $Cay(\Gamma, \{\pm 1\} \times (\pm \{0,1,2\}))$ by
Theorem~\ref{2p-533}. The other $n-7$ $C_n$-factors come
from a $C_n$-factorization of $Cay(\Gamma, \{\pm 1\} \times (\pm \{3,4,\cdots,\frac{n-1}{2}
\}  \backslash \{\pm d\}))$. In fact,
when $n \equiv 3 \pmod {12}$,
there exist 4 $C_n$-factors which form a $C_n$-factorization of
$Cay( \Gamma,  \{ \pm 1\}  \times  (\pm \{ 2j-1,2j \}) )$, $2\leq j \leq \frac{d-1}{2}$ by Theorem~\ref{d1d2}.
There exist 4 $C_n$-factors which form a $C_n$-factorization of
$Cay( \Gamma,  \{ \pm 1\}  \times  (\pm \{ d+2j-1,d+2j \}) )$, $1\leq j \leq \frac{d-1}{4}$ by the same theorem;
when $n \equiv 9 \pmod {12}$,
there exist 4 $C_n$-factors which form a $C_n$-factorization of
$Cay( \Gamma, $ $ \{ \pm 1\}  \times  (\pm \{ 2j-1,2j \}) )$, $2\leq j \leq \frac{d-1}{2}$ by Theorem~\ref{d1d2}.
There exist 4 $C_n$-factors which form a $C_n$-factorization of
$Cay( \Gamma,  \{ \pm 1\}  \times  (\pm \{ d+2j-1,d+2j \}) )$, $1\leq j \leq \frac{d-3}{4}$ by the same theorem.
Further, there exist $2$ $C_n$-factors which form a $C_n$-factorization of
$Cay( \Gamma,  \{ \pm 1\}  \times  \{\pm \frac{n-1}{2}\} )$
by Theorem~\ref{d}.
\qed

\begin{lemma}
\label{Cm[9]4} $(4,5)\in$ \rm{HWP}$(C_m[9];m,9)$ for any odd
integer $m\geq 5$.
\end{lemma}

\noindent {\it Proof:}   Let the vertex set be $\Gamma=Z_{m}\times
Z_9$. Let $C_j=((0,0),(1,b_{j1}),\cdots,(m-1,b_{j,m-1}))$, $1\leq
j\leq 4$, where $ b_{jt}=b_{j,t-2}$, $t\geq 5$ and

\[\begin{array}{lllll}
b_{11}=-b_{21}=3, &  b_{12}=-b_{22}=7, & b_{13}=-b_{23}=3, & b_{14}=-b_{24}=6,  \\
b_{31}=-b_{41}=5, &  b_{32}=-b_{42}=2, & b_{33}=-b_{43}=8, &
b_{34}=-b_{44}=4.
\end{array}
\]

Each $C_j$ will generate a $C_m$-factor by $(-, +1
\pmod{9})$. Thus we get 4 $C_m$-factors which form  a
$C_m$-factorization of $Cay(\Gamma,\{ \pm 1\} \times
(\pm\{3, 4\}))$. The required $5$ $C_9$-factors can be obtained from
a $C_n$-factorization of
$Cay(\Gamma, \{\pm 1\} \times  (\pm \{0,1,2\}))$ by
Theorem~\ref{2p-533}. \qed

\begin{construction}
\label{00}  Let $k\geq3$. If a $k$-{\rm ARCS}$(2kt+1)$ exists, then $Cay(Z_{k}
\times Z_{2kt+1}, \{0\}\times (Z_{2kt+1} \setminus  \{0\}) \cup
\{\pm 1\}  \times \{0\})$ can be decomposed into $kt+1$
$C_k$-factors.
\end{construction}

\noindent {\it Proof:} Let $A_i=\{i\} \times Z_{2kt+1}$, $i \in
Z_{k}$. Put the blocks of a $k$-ARCS$(2kt+1)$ on
each $A_i$. For $0\leq j \leq kt-1$, without loss of generality,
suppose the vertex $(i,j)$ doesn't appear in the $j$th almost parallel class
$P_{ij}$. And $\{(i,j) \ | \ 0\leq j \leq kt-1\}$ appear in the half-parallel class $Q_i$.
Let $S_j=(\bigcup_{i\in Z_k}P_{ij})\bigcup \
((0,j), (1,j),\cdots, (k-1,j))$, $0\leq j \leq kt-1$. It is easy
to check that each $S_j$ is a $C_k$-factor. $T=(\bigcup_{i\in
Z_k}Q_i)\bigcup \  \{ ((0,j), (1,j),\cdots, (k-1,j)) \ | \ kt \leq
j \leq 2kt\}$ is a $C_k$-factor. Namely, we obtain $kt+1$
$C_k$-factors. \qed

In order to give the next construction, we first introduce the following lemma.

\begin{lemma}
\label{17-35} $Cay(Z_{17} \times Z_{35}, \{\pm 1 \}\times
(Z_{35} \setminus  \{0\}) )$ can be decomposed into $28$
$C_{17}$-factors and $6$ $C_{35}$-factors.
\end{lemma}

\noindent {\it Proof:}
Let $A=[a_{ij}]$ be the $28\times 17$ array, where

\[\begin{array}{lllllll}
a_{11}=3, &  a_{12}=-15, & a_{13}=12, &
a_{21}=6, &  a_{22}=12,  & a_{23}=-18,   \\
a_{31}=-12, &  a_{32}=6, & a_{33}=6, &
a_{41}=15,  &  a_{42}=-18, & a_{43}=3,   \\
a_{51}=-18, &  a_{52}=3,   & a_{53}=15, &

a_{61}=2, &  a_{62}=7, & a_{63}=-9,  \\
a_{71}=7,  &  a_{72}=-9, & a_{73}=2, &
a_{81}=-9, &  a_{82}=2,   & a_{83}=7, \\

a_{91}=1, &  a_{92}=10, & a_{93}=-11, &
a_{10,1}=10,  &  a_{10,2}=-11, & a_{10,3}=1,  \\
a_{11,1}=-11, &  a_{11,2}=1,   & a_{11,3}=10, &

a_{12,1}=5, &  a_{12,2}=14, & a_{12,3}=16,  \\
a_{13,1}=14,  &  a_{13,2}=16, & a_{13,3}=5, &
a_{14,1}=16, &  a_{14,2}=5,   & a_{14,3}=14, &
\end{array}
\]

$$a_{ij}=-a_{i-14,j}, \ 15\leq i \leq 28, \ 1\leq j \leq 3, $$
$$a_{i4}=-a_{i5}=a_{i1},  a_{ij}=a_{i,j-2}, \  1\leq i \leq 28, \  6 \leq j \leq 17.$$

Let $b_{ij}=\sum_{l=1}^{j}a_{il}$.
For $1\leq i \leq 28$, each $C_i=(0_0, 1_{b_{i1}}, 2_{b_{i2}}, \cdots, 16_{b_{i,16}})$ can get a $C_{17}$-factor by $(-, +1 \pmod {35})$.
It is easy to check that these $28$ $C_{17}$-factors are a $C_{17}$-factorization
of $Cay(Z_{17} \times Z_{35}, \{\pm 1 \}\times
(Z_{35} \setminus (\pm \{0,4,8,13\}) ))$.

For $C_{35}$-factors, four of which come from a $C_{35}$-factorization
of $Cay(Z_{17} \times Z_{35}, \{\pm 1 \}\times (\pm  \{4,8\}) )$ by Theorem~\ref{d1d2}.
The other 2 $C_{35}$-factors come from a $C_{35}$-factorization
of $Cay(Z_{17} \times Z_{35}, \{\pm 1 \}\times ( \{\pm 13\}) )$ by Theorem~\ref{d}.
\qed

\begin{construction}
\label{2l} Let $k$ be an odd integer. Then $Cay(Z_{k} \times
Z_{2kt+1}, \{\pm 1 \}\times (Z_{2kt+1} \setminus  \{0\}) )$ can be
decomposed into $2l$ $C_k$-factors and $2kt-2l$
$C_{2kt+1}$-factors for $l \notin \{1,2,kt-1,kt\}$.
\end{construction}

\noindent {\it Proof:} By Lemmas 3.6, 3.11 and 3.13 in \cite{BDT}
we know that $Cay(Z_{k} \times Z_{2kt+1}, \{\pm 1 \}\times
Z_{2kt+1} )$ can be decomposed into $2l$ $C_k$-factors and
$2kt-2l+1$ $C_{2kt+1}$-factors for any $l \notin \{1,2,kt-1,kt\}$.
In these constructions, we can always find five
$C_{2kt+1}$-factors which form a $C_{2kt+1}$-factorization of
$Cay(Z_{k} \times Z_{2kt+1}, \{\pm 1 \}\times(\pm   \{0,d,2d\}) )$
such that $(d,2kt+1)=1$ except for $(k,t,l)=(17,1,14)$. By
Theorem~\ref{d1d2} we know that $Cay(Z_{k} \times Z_{2kt+1}, \{\pm
1 \}\times(\pm   \{d,2d\}) )$ can be partitioned into $4$
$C_{2kt+1}$-factors. Thus $Cay(Z_{k} \times Z_{2kt+1}, \{\pm 1
\}\times (Z_{2kt+1} \setminus  \{0\}) )$ with
$(k,t,l)\not=(17,1,14)$ can be decomposed into $2l$ $C_k$-factors
and $2kt-2l$ $C_{2kt+1}$-factors for $l \notin \{1,2,kt-1,kt\}$.
Further, $Cay(Z_{17} \times Z_{35}, \{\pm 1 \}\times
(Z_{35} \setminus  \{0\}) )$ can be decomposed into $28$
$C_{17}$-factors and $6$ $C_{35}$-factors by Lemma~\ref{17-35}.
\qed

\begin{construction}
\label{2ku} If a $k$-{\rm CF}$(2^u)$ exists, then $Cay(Z_{k}
\times Z_{2u}, \{0\}\times (Z_{2u} \setminus  \{0\}) \cup \{\pm
1\} \times \{0\})$ can be decomposed into $u$ $C_k$-factors and a
$1$-factor.
\end{construction}

\noindent {\it Proof:}
Let $A_i=\{i\} \times Z_{2u}$, $i \in Z_{k}$.
The $1$-factor is $\{ (i,2j), (i,2j+1)   \ | \ i \in Z_{k}, 0\leq j \leq u-1 \}$
For each $A_i$ , put the blocks of a $k$-{\rm CF}$(2^u)$ on
the vertex set $A_i$, the parts of the cycle frame are
$G_{ij}=\{ (i,2j), (i,2j+1) \}$, $0\leq j \leq u-1$. Denote the holey 2-factor of $G_{ij}$ by $P_{ij}$.
 Let $S_j=(\bigcup_{i\in Z_k}P_{ij})\bigcup \
\{ ((0,2j), (1,2j),\cdots, (k-1,2j)), ((0,2j+1), (1,2j+1),\cdots, (k-1,2j+1))\}$, $0\leq j \leq u-1$. It is easy
to check that each $S_j$ is a $C_k$-factor. That is to say, we get $u$
$C_k$-factors. \qed

 \section{Main results}

In order to prove the main results, we need the following six
direct constructions.

\begin{lemma}
\label{C33} $(6,5)\in$ \rm{HWP}$(C_3[11];3,11)$.
\end{lemma}

\noindent {\it Proof:}   Let the vertex set be $(Z_{6}\times
Z_5)\cup \{ a,b,c \}$, and the three parts of $C_3[11]$ be $\{ a,
\{2, 5\} \times Z_5  \}$, $\{ b, \{0, 3\} \times Z_5 \}$ and $\{ c, \{1, 4\} \times Z_5 \}$.

For the required 6 $C_3$-factors, five of which will be generated
from an initial $C_3$-factor $P$ by $(-, +1 \pmod{5})$. The last
$C_3$-factor is $\{( a, b, c), (0_i, 1_{2+i}, 2_{4+i}),( 3_{3+i}, 4_{2+i}, 5_{3+i}) \ | \ i \in Z_5 \}$.
5 $C_{11}$-factors will be generated from an initial
$C_{11}$-factor $Q$ by $(-, +1 \pmod{5})$. The cycles of $P$ and
$Q$ are listed below.

{\footnotesize \vspace{5pt}\noindent
\begin{tabular}{lllllllllll}
$P$ &$( a, 0_0, 4_4)$ &$( b, 2_3, 1_0)$ &$( c, 2_4, 3_1)$ &$( 1_1,
3_3, 5_0)$ &$( 2_2, 4_0, 0_2)$
&$( 0_1, 5_2, 1_4)$\\
&$( 1_2, 2_1, 3_2)$ &$( 3_4, 2_0, 4_2)$ &$( 5_1, 0_3, 4_3)$ &$(
1_3, 0_4, 5_4)$
&$( 3_0, 4_1, 5_3)$\\
\end{tabular}}

{\footnotesize \noindent
\begin{tabular}{lllllllll}
$Q$ &$(a, 1_1, 5_1, 3_0, 1_4, 0_4, 4_1, 2_2, 3_2, 1_3, 3_1)$
&$(b, 5_0, 0_1, 1_2, 2_3, 4_3, 0_2, 2_4, 0_3, 2_1, 4_2)$\\
&$(c, 0_0, 5_2, 1_0, 2_0, 3_3, 4_0, 5_4, 4_4, 3_4, 5_3)$
\end{tabular}}

 \qed

\begin{lemma}
\label{C39} $(8,5)\in$ \rm{HWP}$(C_3[13];3,13)$.
\end{lemma}

\noindent {\it Proof:}   Let the vertex set be $\Gamma=Z_{13}\times Z_3$,
and the three parts of $C_3[13]$ be $Z_{13}\times \{i\}$, $i\in
Z_3$.

For $C_3$-factors, five of which come from a $C_3$-factorization of $Cay(\Gamma, \pm
\{0,1,2\}\times \{\pm 1\})$ by Lemma~\ref{2p-533}.
 The 3-cycle $( 0_0, 4_1, 7_2)$ can generate a
$C_3$-factor $F$ by $(+1 \pmod {13}, -)$, then the last 3
$C_3$-factors can be obtained from $F$ by $(-,+1 \pmod {3})$.
5 $C_{13}$-factors will be generated from 5 13-cycles in $Q$ by
$(-, +1 \pmod{3})$ since the 13 elements in the first coordinate
of each 13-cycle are different modular 13. The cycles of
$Q$ are listed below.

{\footnotesize \vspace{5pt}\noindent
\begin{tabular}{lllllllllll}
$Q$ &$(0_0, 5_2, 10_1, 2_2, 7_1, 12_0, 4_1, 9_0, 1_1, 6_0, 11_2,
3_0, 8_2)$
&$(0_0, 7_1,  1_2, 4_1, 11_2, 3_1, 6_0, 10_2, 2_1, 8_0, 12_2, 9_0, 5_1)$\\
&$(0_0, 10_1, 3_0, 7_2, 12_0, 2_2, 9_0, 1_2, 6_0, 11_1, 5_2, 8_1,
4_2)$
&$(0_0, 9_1,  3_2, 8_0, 11_2, 1_1, 10_2, 4_0, 7_2, 2_1, 5_0, 12_1, 6_2)$\\
&$(0_0, 3_2, 12_0, 4_2, 9_0,  6_1, 2_2, 11_0, 7_1, 10_0, 5_2, 1_0,
8_1)$
\end{tabular}}

 \qed

\begin{lemma}
\label{C45} $(8,7)\in$ \rm{HWP}$(C_3[15];3,15)$.
\end{lemma}

\noindent {\it Proof:}   Let the vertex set be $Z_{45}$, and the
three parts of $C_3[15]$ be $\{3i+j \ | \ i \in Z_{15}\}$,
$j\in Z_3$. Since all 3 elements of each cycle in $P$ are
different modular 3, so all 8 $C_{3}$-factors will be generated
from 8 3-cycles in $P$ by $+3 \pmod{45}$. And the required 7
$C_{15}$-factors will be generated from 7 15-cycles in $Q$ by $+15
\pmod{45}$. The cycles of $P$ and $Q$ are listed below.

{\footnotesize \vspace{5pt}\noindent
\begin{tabular}{lllllllllll}
$P$ &$( 0, 1, 2 )$ &$( 0, 4, 8 )$ &$( 0, 5, 7 )$ &$( 0, 10, 17 )$
&$( 0, 11, 16 )$ &$( 0, 13, 23 )$ &$( 0, 14, 22 )$ &$( 0, 19, 32
)$
\end{tabular}}

{\footnotesize  \noindent
\begin{tabular}{lllllllll}
$Q$ &$( 0, 34, 6, 41, 3,  43, 9,  44, 10, 27, 2,  16, 5,  22, 38
)$
&$( 0, 20, 1, 3,  2,  19, 21, 14, 43, 26, 40, 9,  38, 42, 37 )$\\
&$( 0, 28, 2, 12, 1,  18, 7,  21, 5,  9,  8,  19, 44, 25, 41 )$
&$( 0, 29, 1, 17, 28, 3,  23, 4,  20, 39, 22, 42, 26, 6,  40 )$\\
&$( 0, 31, 6, 32, 3,  29, 9,  34, 11, 37, 23, 42, 40, 20, 43 )$
&$( 0, 25, 2, 9,  1,  6,  5,  19, 33, 23, 43, 12, 11, 22, 44 )$\\
&$( 0, 26, 1, 23, 40, 3,  44, 13, 21, 17, 37, 39, 34, 42, 35 )$
\end{tabular}}

 \qed

Combining Theorem~\ref{Cn-m}, Lemmas~\ref{C33}-\ref{C45}, we
have the following conclusion.

\begin{lemma}
\label{Cn-m2}  If $n$, $m$ and $\alpha$ are odd integers with $n
\geq m \geq 3$, $0 \leq \alpha \leq n$, then $(\alpha,\beta) \in $
{\rm HWP}$(C_m[n]; m, n)$, if and only if $\beta = n-\alpha$,
except possibly when $\alpha=2,4$, $\beta=1,3$.
\end{lemma}

\begin{lemma}
\label{K33} $(6,10)\in$ \rm{HWP}$(33;3,11)$.
\end{lemma}

\noindent {\it Proof:}   Let the vertex set be $(Z_{6}\times
Z_5)\cup \{ a,b,c \}$. For the required 6 $C_3$-factors, five of which
will be generated from an initial $C_3$-factor $P$ by $(-, +1
\pmod{5})$. The last $C_3$-factor is
$\{( a, b, c), ( (3j)_i, (3j+1)_{3+i}, (3j+2)_i) \ | \ i \in Z_5, j=0,1\}$.
10 $C_{11}$-factors will be generated from an
initial $C_{11}$-factor $Q$ by $( +3 \pmod{6}, +1 \pmod{5})$.
$P$ and $Q$ are listed below.

{\footnotesize \vspace{5pt}\noindent
\begin{tabular}{lllllllllll}
$P$ &$( a, 1_1, 4_1)$ &$( b, 2_2, 5_2)$ &$( c, 0_0, 3_0)$ &$( 3_3,
5_0, 1_4)$ &$( 4_4, 0_3, 2_0)$
&$( 0_1, 2_4, 3_2)$\\
&$( 1_2, 2_3, 2_1)$ &$( 3_4, 4_0, 4_3)$ &$( 5_1, 4_2, 5_3)$ &$(
0_2, 3_1, 5_4)$
&$( 1_3, 0_4, 1_0)$\\
\end{tabular}}

{\footnotesize  \noindent
\begin{tabular}{lllllllll}
$Q$ &$( a, 2_3, 1_3, 0_3, 5_1, 4_3, 3_1, 2_0, 3_4, 3_0, 3_2)$
&$( b, 0_2, 2_1, 5_4, 0_4, 4_1, 4_2, 1_0, 5_2, 5_3, 1_4)$\\
&$( c, 1_1, 4_0, 0_0, 3_3, 1_2, 5_0, 2_4, 4_4, 0_1, 2_2)$
\end{tabular}}

 \qed

 \begin{lemma}
\label{K35} $(9,8)\in$ \rm{HWP}$(35;5,7)$.
\end{lemma}

\noindent {\it Proof:}   Let the vertex set be $\Gamma=Z_{5}\times
Z_7$. For the required 9 $C_5$-factors, seven of which will be
generated an initial $C_5$-factor $P$ by $(-, +1 \pmod{7})$. The
last 2 $C_5$-factors will be obtained from two cycles
 $\{( 0_0, 4_2, 2_2, 1_4, 3_4), ( 0_0, 1_4, 3_3, 2_3, 4_1)\}$
 by $(-, +1 \pmod{7})$ since the 5 elements in the first coordinate
 of each cycle are different modular 5.
For $C_7$-factors, seven of which will be generated an initial
$C_7$-factor $Q$ by $(-, +1 \pmod{7})$. The last $C_7$-factor is
$Cay(\Gamma, \{0\}\times \{\pm 2\})$. The cycles of $P$ and $Q$
are listed below.

{\footnotesize \vspace{5pt}\noindent
\begin{tabular}{lllllllllll}
$P$ &$( 0_0, 1_1, 2_2, 3_3, 4_4)$ &$( 0_5, 2_0, 4_2, 1_6, 3_1)$
&$( 0_3, 1_2, 2_5, 0_1, 2_1)$
&$( 1_4, 2_3, 3_6, 4_5, 0_6)$\\
&$( 4_0, 3_4, 0_4, 3_5, 2_6)$ &$( 1_0, 4_1, 2_4, 3_2, 4_6)$ &$(
4_3, 1_5, 0_2, 3_0, 1_3)$
\end{tabular}}

{\footnotesize  \noindent
\begin{tabular}{lllllllll}
$Q$ &$( 0_0, 0_3, 2_2, 2_5, 1_1, 1_4, 3_2)$ &$( 3_3, 3_6, 1_5,
0_5, 2_1, 4_4, 1_0)$
&$( 1_6, 2_6, 2_0, 0_6, 4_6, 4_5, 4_1)$\\
&$( 3_1, 4_3, 2_4, 1_2, 1_3, 0_1, 0_2)$ &$( 4_2, 3_4, 3_5, 2_3,
3_0, 4_0, 0_4)$
\end{tabular}}

 \qed

\begin{lemma}
\label{K39} $(8,11)\in$ \rm{HWP}$(39;3,13)$.
\end{lemma}

\noindent {\it Proof:}   Let the vertex set be
$\Gamma=Z_{13}\times Z_3$. For $C_3$-factors, five of which will be
generated from a $C_3$-factorization of $Cay(\Gamma, \pm
\{0,1,2\}\times \{\pm 1\})$ by Lemma~\ref{2p-533}. The 3-cycle $(
0_0, 4_1, 7_2)$ can generate a $C_3$-factor $F$ by $(+1 \pmod
{13}, -)$, then the last 3 $C_3$-factors can be obtained from $F$
by $(-,+1 \pmod {3})$. 11 $C_{13}$-factors will be generated
from 11 13-cycles in $Q$ by $(-,+1 \pmod {3})$. The cycles of $Q$
are listed below.

{\footnotesize \vspace{5pt}\noindent
\begin{tabular}{lllllllll}
$Q$ &$( 0_0, 4_2, 2_2, 6_1, 1_1, 3_1,11_0, 7_1, 5_1, 8_0,12_2,
9_0,10_0)$
&$( 0_0, 5_0, 3_0, 7_2, 2_2,10_2, 4_0, 1_1, 6_2,11_0,12_0, 9_0, 8_0)$\\
&$( 0_0, 9_1, 2_1, 5_0, 4_0, 8_0, 6_0,10_0, 1_2,11_2, 7_2,12_2,
3_2)$
&$( 0_0, 3_0, 6_0, 1_1, 4_1, 9_2, 5_2,10_2,12_2, 7_1, 8_1, 2_2,11_0)$\\
&$( 0_0, 5_2,10_1, 1_1, 9_1, 3_1, 6_0,11_2, 8_0, 2_0, 7_1,12_0,
4_0)$
&$( 0_0, 7_1,10_1, 8_1,11_1, 3_0, 4_0, 9_0, 1_1, 2_1, 5_1,12_1, 6_2)$\\
&$( 0_0,10_1, 2_0, 3_0,11_0, 5_1, 8_1, 4_2,12_1, 6_1, 9_1, 1_0,
7_0)$
&$( 0_0,12_0, 2_0, 9_1, 4_2,11_0,10_0, 6_1, 7_1, 3_1, 8_2, 1_1, 5_1)$\\
&$( 0_0, 2_0, 6_0, 5_0,10_1, 7_2, 4_0,12_1, 8_1, 3_2, 9_1,11_1,
1_0)$
&$( 0_0, 8_2, 1_2,12_2, 3_1,10_1, 2_2, 7_1, 4_1, 6_1,11_1, 5_1, 9_0)$\\
&$( 0_0, 6_0, 9_2, 7_2, 1_0, 5_2,12_0, 2_2,11_2, 4_2,10_2, 3_1,
8_1)$

\end{tabular}}

 \qed

\begin{lemma}
\label{63}   $(\alpha,\beta)\in$ \rm{HWP}$(9u;u,9)$ for $u=5,7$, $\beta=9,11$ and $\alpha+\beta=\frac{9u-1}{2}$.
\end{lemma}

\noindent {\it Proof:} There exist an HW$(K_u[1];u,9;
\frac{u-1}{2},0)$ (from Theorem~\ref{Kv}), an HW$(C_u[9];u,9;9,0)$
(from Theorem~\ref{Cn-m}) and an HW$(C_u[9];u,9;\alpha,9-\alpha)$
for $\alpha=2,4$ (from Lemmas~\ref{Cm[n]2} and~\ref{Cm[9]4}).
Apply Construction~\ref{L351} with $s=9$ and $t_i\in\{9,\alpha\}$
to obtain an HW$(K_u[9];u,9; \frac{9(u-3)}{2}+\alpha,9-\alpha)$.
Applying Construction~\ref{C-RGDD} with an HW$(9;u,9;0,4)$ from
Theorem~\ref{Kv}, we can get an HW$(9u;u,9;
\frac{9(u-3)}{2}+\alpha,13-\alpha)$ for $\alpha=2,4$. \qed

\begin{lemma}
\label{39t}  $(\frac{39t-11}{2}, 5)\in$ \rm{HWP}$ (39t; 3, 13)$  for any odd $ t>1$.
\end{lemma}

\noindent {\it Proof:} Let $\alpha=\frac{3(t-1)}{2}$. Start with
an HW$(K_t[3];3,13;\alpha,0)$ from Theorem~\ref{Kv}, an
HW$(C_3[13];$ $3,13;13,0)$ from Theorem~\ref{Cn-m} and an
HW$(C_3[13];3,13;8,5)$ from Lemma~\ref{C39}. Then apply
Construction~\ref{L351} with $s=13$ and $t_i\in\{8,13\}$ to get an
HW$(K_t[39];3,13;13\alpha-5,5)$. Apply Construction~\ref{C-RGDD}
with an HW$(39;3,13;19,0)$ from Theorem~\ref{Kv} to get the
conclusion. \qed

\begin{lemma}
\label{63t}  $(\alpha, \beta)\in$ \rm{HWP}$ (9tu; u, 9)$  for any odd $ t>1$, $u=\beta=5,7$ and $\alpha+\beta=\frac{9tu-1}{2}$.
\end{lemma}

\noindent {\it Proof:}
There exist an
HW$(K_t[u];u,9;\alpha',0)$, $\alpha'=\frac{u(t-1)}{2}$ (from Theorem~\ref{Kv}), an HW$(C_u[9];u,$ $9;9,0)$ (from Theorem~\ref{Cn-m}) and an
HW$(C_u[9];u,9;\alpha,9-\alpha)$, $\alpha=2,4$ (from Lemmas~\ref{Cm[n]2} and~\ref{Cm[9]4}). Apply
Construction~\ref{L351} with $s=9$ and $t_i\in\{9,\alpha\}$ to get an
HW$(K_t[9u];u,9;9(\alpha'-1)+\alpha,9-\alpha)$.
Applying Construction~\ref{C-RGDD} with an
HW$(9u;u,9;\frac{9u-1}{2},0)$ from Theorem~\ref{Kv}, the assertion follows.
\qed

\noindent {\bf  Proof of Theorem~\ref{odd2}:} Combining
Theorem~\ref{odd} and Lemmas~\ref{63}-\ref{63t}, we can get the
conclusion. \qed

\noindent {\bf  Proof of Theorem~\ref{ARCS}:}
By Theorem~\ref{Kv}
the complete graph $K_k$ with vertex set $Z_k$ can be decomposed
into $r=\frac{k-1}{2}$ cycles $B_1,B_2,\ldots, B_r$. Without lose
of generality, we may suppose that $B_r=(0,1,2,\ldots,k-1)$. Give
each vertex weight $2kt+1$. Let $\Gamma=Z_{k} \times Z_{2kt+1}$.

When $k=3$ and $t\geq 3$, we have $r=1$.
For the cycle $B_r$, $Cay(\Gamma, (\{0\}\times (Z_{6t+1}
\setminus \{0\})) \cup  (\{ \pm 1\} \times \{0\}))$ can be
decomposed into $3t+1$ $C_3$-factors by Construction~\ref{00},
where a $3$-{\rm ARCS}$(6t+1)$ for $t\geq3$ exists by Theorem~\ref{3-14}.
By Construction~\ref{2l} $Cay(\Gamma, \{\pm 1
\}\times (Z_{6t+1} \setminus  \{0\}) )$ can be decomposed into
$2l$ $C_3$-factors and $6t-2l$ $C_{6t+1}$-factors, $ 2l \in \{0, 6, 8, 10, \cdots, 6t-4\}$.
There are altogether $\alpha=2l+3t+1$ $C_3$-factors
and $\beta=6t-2l$ $C_{6t+1}$-factors.
It is obvious that $\beta \in [4, 5£¬ \cdots, 6t-6] \cup \{6t\}$ and $\beta$ is even.

Next, we consider the case $k\geq 5$. Then we have $r\geq 2$. For each cycle $B_i$, $1 \leq i \leq r-1$, the graph $C_k[2kt+1]$
can be decomposed into $2kt+1-\beta_i$ $C_k$-factors and $\beta_i$
$C_{2kt+1}$-factors by Theorem~\ref{Cn-m}, where $0\le \beta_i\le
2kt+1$, $\beta_i \not\in \{1,3,2kt-3,2kt-1\}$.
For the cycle $B_r$, $Cay(\Gamma, (\{0\}\times (Z_{2kt+1}
\setminus \{0\})) \cup  (\{ \pm 1\} \times \{0\}))$ can be
decomposed into $kt+1$ $C_k$-factors by Construction~\ref{00},
where a $k$-{\rm ARCS}$(2kt+1)$ for $t\geq 1$ ($t\neq 2$ when $k \geq 11$) exists by Theorem~\ref{3-14} and
Lemma~\ref{2kt+1}.  By Construction~\ref{2l} $Cay(\Gamma, \{\pm 1
\}\times (Z_{2kt+1} \setminus  \{0\}) )$ can be decomposed into
$2l$ $C_k$-factors and $2kt-2l$ $C_{2kt+1}$-factors, $0\le 2l\le
2kt+1$, $l \not\in \{1,2,kt-1,kt\}$.
Here we have obtained
$2l+kt+1$ $C_k$-factors and $2kt-2l$ $C_{2kt+1}$-factors.
Finally, we get $\alpha=(2l+kt+1)+ \sum_{i=1}^{r-1}(2kt+1-\beta_i)$ $C_k$-factors
and $\beta=(2kt-2l)+ \sum_{i=1}^{r-1}\beta_i$ $C_{2kt+1}$-factors, where $0\le \beta_i, 2l \le
2kt+1$, $\beta_i \not\in \{1,3,2kt-3,2kt-1\}$ and  $l \not\in \{1,2,kt-1,kt\}$.
When $k=5$, then $r=2$ and $\beta =(10t-2l)+\beta_1$, where $2l \in \{0, 6, 8, 10, \cdots, 10t-4\}$
and $\beta_1=[0, 1, \cdots, 10t+1]\backslash \{1, 3, 10t-3, 10t-1\}$. It is easy to check that $\beta \in J \setminus \{ 20t-3, 20t-1 \}$.
When $k\geq 7$, then $\beta=(2kt-2l)+ \sum_{i=1}^{r-1}\beta_i$, where $2l \in \{0, 6, 8, 10, \cdots, 2kt-4\}$
and $\beta_i \in [0, 1, \cdots, 2kt+1]\backslash \{1,3,2kt-3,2kt-1\}$. We can check that $\beta \in J $.
The proof is complete.
\qed

\noindent {\bf  Proof of Theorem~\ref{4k-4kt}:} Let $v=4ktu$,
$u\geq 1$. For $u=1$,  the conclusion comes from Theorem 1.6 in \cite{LFS}.
For $u\geq 2$, start with an HW$(K_u[4k];4k,4kt;2k(u-1),0)$, an
HW$(C_{4k}[t];4k,4kt;t,$ $0)$ and an HW$(C_{4k}[t];4k,4kt;0,t)$ from
Theorems~\ref{Kv},~\ref{Cm} and~\ref{Cmn}, respectively.
 And apply Construction~\ref{L351} with $s=t$ and
$t_i\in\{0,t\}$ to get an
HW$(K_u[4kt];4k,4kt;\sum_{i=1}^{2k(u-1)}t_i,2kt(u-1)-\sum_{i=1}^{2k(u-1)}t_i)$.
Further, applying Construction~\ref{C-RGDD} with an
HW$(4kt;4k,4kt;\alpha',2kt-1-\alpha')$, $0\leq \alpha'\leq2kt-1$, from
Theorem 1.6 in \cite{LFS},
 we can obtain an
HW$(4ktu;4k,4kt;\alpha'+\sum_{i=1}^{2k(u-1)}t_i,(2ktu-1)-(\alpha'+\sum_{i=1}^{2k(u-1)}t_i))$.
It's easy to check that $\alpha'+\sum_{i=1}^{2k(u-1)}t_i$ can
cover all the integers from $0$ to $2ktu-1$.  \qed

 \section{Concluding remarks}

Combining Theorems~\ref{odd2} and \ref{ARCS}, we have the following open problem for $m=k$ and $n=2kt+1$.

\begin{problem}
\label{Q1}
Find a solution \rm{HW}$ (k(2kt+1); k, 2kt+1; \alpha,\beta)$ in the following cases:

$(i)$ $t=1$

$\ \ \ 1)$ $k=5$: $\beta = \{1, 2, 3\}$;

$\ \ \ 2)$ $k=7$: $\beta \in \{ 1, 2, 3, 5\}$;

$\ \ \ 3)$ $k\geq 9$: $\beta \in \{ 1, 2, 3, 5, 7\}$.

$(ii)$ $t=2$

$\ \ \ 1)$ $k=3$ or $k\geq 11$: $\beta \in [1,2,\cdots, 2k-1]\cup \{2k+1, 2k+3\}$;

$\ \ \ 2)$ $k=5, 7, 9 $: $\beta \in \{ 1, 2, 3, 5, 7\}$.

$(iii)$ $t\geq 3$

$\ \ \ 1)$ $k=3$:

$\ \ \ \  \bullet$ $t$ is odd: $\beta= \{1, 3, 5, \cdots, 3t-2\} \cup \{2, 9t-3, 9t-1\}$

$\ \ \ \  \bullet$ $t$ is even: $\beta= \{1, 3, 5, \cdots, 3t+3\} \cup \{2, 9t-3, 9t-1\}$;

$\ \ \ 2)$ $k\geq 5$: $\beta \in \{ 1, 2, 3, 5, 7\}$.
\end{problem}


\begin{thebibliography}{Z}
\baselineskip 11pt

\bibitem{ABBE}
 P. Adams, E. J. Billington, D. E. Bryant, S.I. El-Zanati, On the Hamilton-Waterloo problem,
 {\it Graphs Combin.} {\bf 18} $(2002)$, 31-51.


\bibitem{ABHL}
 P. Adams, E. J. Billington, D. G. Hoffman and C. C. Lindner, {The
generalized almost resolvable cycle system problem},
{\it Combinatorica} {\bf 30} (2010), 617-625.


 \bibitem{AG}
 B. Alspach, H. Gavlas, Cycle decompositions of $K_n$ and $K_n- I$,
 {\it J. Combin. Theory Ser. B} {\bf 81} $(2001)$, 77-99.

\bibitem{AH}
A. Assaf, A. Hartman, Resolvable group divisible designs with
block size $3$, {\it Discrete Math. }{\bf 77} $(1989)$,  5-20.


 \bibitem{AKK}
 J. Asplund, D. Kamin, M. Keranen, A. Pastine, S. Ozkan, On the  Hamilton-Waterloo problem with triangle factors and $C_{3x}$-factors,
$(2015)$, arXiv:$1510.04607$ [math.CO].


\bibitem{ASSW}
B. Alspach,  P. J. Schellenberg,  D. R. Stinson, D. Wagner, The
Oberwolfach problem and factors of uniform odd length cycles, {\it
J. Combin.  Theory Ser.  A} {\bf 52} $(1989)$,  20-43.


\bibitem{BCDT}
M. Buratti, H. Cao, D. Dai and  T. Traetta, A complete solution to
the existence of $(k,\lambda)$-cycle frames of type $g^u$,
preprint.

 \bibitem{BDT}
 A. Burgess, P. Danziger, T. Traetta, On the Hamilton-Waterloo problem with odd orders,
 $(2015)$, arXiv:$1510.07079$ [math.CO].


 \bibitem{BH}
J. Burling and K. Heinrich, Near $2$-factorizations of $2K_{n}$:
cycles of even length, {\it Graphs Combin.} {\bf 5}(1989),
213-221.


 \bibitem{BHL}
 E. J. Billington, D. G. Hoffman, C. C. Lindner, M. Meszka, Almost resolvable minimum coverings of complete graphs with 4-cycles,
 {\it  Australas. J. Combin.} {\bf 50} (2011), 73-85.



\bibitem{CNT}
H. Cao,  M. Niu, C. Tang, On the existence of cycle frames and
almost resolvable cycle systems, {\it Discrete Math.} {\bf
311} (2011), 2220-2232.



 \bibitem{DQS}
 P. Danziger, G. Quattrocchi, B. Stevens, The Hamilton-Waterloo problem for cycle sizes $3$ and $4$,
 {\it J. Combin. Des.} {\bf 17} $(2009)$, 342-352.

 \bibitem{DL}
 J. H. Dinitz, A. C. H. Ling, The Hamilton-Waterloo problem with triangle-factors and Hamilton cycles: The case $n\equiv 3\pmod{18}$,
 {\it J. Combin. Math. Combin. Comput. } {\bf 70} $(2009)$, 143-147.

 \bibitem{DL2}
 J. H. Dinitz, A. C. H. Ling, The Hamilton-Waterloo problem: The case of triangle-factors and one Hamilton cycle,
 {\it J. Combin. Des.} {\bf 17} $(2009)$, 160-176.

  \bibitem{DLM}
  I. J. Dejter, C. C. Lindner, M. Meszka, C. A. Rodger, Corrigendum/addendum to: almost resolvable 4-cycle systems,
 {\it J. Combin. Math. Combin. Comput. } {\bf 66} $(2008)$, 297-298.

  \bibitem{DLR}
  I. J. Dejter, C. C. Lindner, C. A. Rodger, M. Meszka, Almost resolvable 4-cycle systems,
 {\it J. Combin. Math. Combin. Comput. } {\bf 63} $(2007)$, 173-181.


\bibitem{ER}
A. Erzurumluo$\check{\hbox{g}}$lu and  C. A. Rodger, Fair holey
hamiltonian decompositions of complete multipartite graphs and
long cycle frames, {\it Discrete Math.} {\bf 338}(2015),
1173-1177.


\bibitem{HLR}
K. Heinrich, C. C. Lindner and C. A. Rodger, Almost resolvable
decompositions of $2K_{n}$ into cycles of odd length, {\it J.
Combin. Theroy Ser. A } {\bf 49}(1988), 218-232.


 \bibitem{HNR}
 P. Horak, R. Nedela, A. Rosa, The Hamilton-Waterloo problem: The case of Hamilton cycles and triangle-factors,
 {\it Discrete Math.} {\bf 284} $(2004)$, 181-188.

  \bibitem{HS}
 D. G. Hoffman, P. J. Schellenberg, The existence of $C_{k}$-factorizations of $K_{2n}-F$,
 {\it Discrete Math.} {\bf 97} $(1991)$, 243-250.


  \bibitem{K}
 D. C. Kamin, Hamilton-Waterloo problem with triangle and $C_9$-factors,
 Master's Thesis Michigan Technological University, 2011.
 https://digitalcommons/etds/207


\bibitem{KO}
M. Keranen, S. $\ddot{O}$zkan,  The Hamilton-Waterloo problem with 4-cycles and a singer factor of $n$-cycles,
 {\bf 29} (2013), 1827-1837.


\bibitem{JLA}
J. Liu,  A generalization of the Oberwolfach problem and
$C_{t}$-factorizations of complete equipartite graphs, {\it  J.
Combin. Des.} {\bf 8} (2000), 42-49.



  \bibitem{L}
 J. Liu, The equipartite Oberwolfach problem with uniform tables,
 {\it J. Combin. Theory Ser.}  {\it A} {\bf 101} $(2003)$, 20-34.

 \bibitem{LD}
J. Liu,  D. R. Lick,  On $\lambda$-fold equipartite Oberwolfach
problem with uniform table sizes, {\it Ann. Comb.}  {\bf 7} (2003),
315-323.


 \bibitem{LF}
 H. C. Lei, H. L. Fu, The Hamilton-Waterloo problem for triangle-factors and heptagon-factors,
 {\it Graphs Combin.}  {\bf 32} $(2016)$, 271-278.


  \bibitem{LFS}
 H. C. Lei, H. L. Fu, H. Shen, The Hamilton-Waterloo problem for Hamilton cycles and $C_{4k}$-factors,
 {\it Ars. Combin.}  {\bf 100} $(2011)$, 341-347.


\bibitem{LMR}
C. C. Lindner, M. Meszka, A. Rosa, Almost resolvable cycle systems-an alogue of Hanani triple systems,
 {\it J. Combin. Des.}  {\bf 17} $(2009)$, 404-410.


 \bibitem{LS}
 H. C. Lei, H. Shen, The Hamilton-Waterloo problem for Hamilton cycles and triangle-factors,
 {\it J. Combin. Des.} {\bf 20} $(2012)$, 305-316.


 \bibitem{MS}
A. Muthusamy and A. Shanmuga Vadivu, Cycle frames of complete
multipartite multigraphs - III, {\it J. Combin. Des.} {\bf
22}(2014), 473-487.

 \bibitem{NC}
M. X. Niu, H. T. Cao, More results on cycle frames and almost resolvable cycle systems, {\it Discrete Math.} {\bf
312} (2012), 3392-3405.

\bibitem{OO}
 U. Odaba\c{s}{\scriptsize I}, S. $\ddot{{\rm O}}$zkan, The Hamilton-Waterloo problem with $C_{4}$ and $C_{m}$ factors,
 {\it Discrete Math.} {\bf 339} $(2016)$, 263-269.



\bibitem{PWL}
W. L. Piotrowski, The solution of the bipartite analogue of the
Oberwolfach problem, {\it Discrete Math.}  {\bf 97} (1991),
339-356.


\bibitem{R}
R. Rees,  Two new direct product-type constructions for resolvable
group-divisible designs, {\it J. Combin. Des.}  {\bf 1} (1993),
15-26.


\bibitem{S}
D. R. Stinson, Frames for Kirkman triple systems, {\it Discrete
Math.} {\bf 65}(1987),  289-300.


\bibitem{SM}
M. $\check{S}$ajna, Cycle decompositions: complete graphs and fixed length cycles,
{\it J. Combin. Des.}  {\bf 10} (2002), 27-78.


\bibitem{VSS}
S. A. Vanstone, D. R. Stinson, P. J. Schellenberg, A. Rosa, R. Rees, C. J. Colbourn, M. W. Carter, J. E. Carter,
 Hanani triple systems, {\it Israel J. Math.}  {\bf 83} (1993), 305-319.

\bibitem{W}
D. West,  Introduction to Graph Theory, 2nd Edition, Prentice
Hall, 2001.


\bibitem{WCC}
L. Wang, F. Chen, H. T. Cao, The Hamilton-Waterloo problem for $C_3$-factors and $C_n$-factors,
preprint.


\bibitem{ZC}
L. Zhang, More results on cycle frames, Master Thesis, Nanjing
Normal University, 2013.

\end{thebibliography}
\end{document}